\documentclass[11pt]{article}
\usepackage{cite}
\usepackage{mathrsfs}
\usepackage{amsfonts}
\usepackage{amsmath}
\usepackage{amsfonts,amssymb,color}
\usepackage{dsfont}
\usepackage{curves}
\usepackage{mathrsfs}
\usepackage{pifont}
\usepackage{amssymb}
\usepackage{latexsym,amsmath,amssymb,amsfonts,epsfig,graphicx,cite,psfrag}
\usepackage{eepic,color,colordvi,amscd}

\newtheorem{theorem}{Theorem}[section]

\newtheorem{lemma}{Lemma}[section]
\newtheorem{corollary}{Corollary}[section]

\newtheorem{claim}{Claim}[section]

\newcommand{\qed}{\hfill\rule{0.5em}{0.809em}}

\def\emptyset{\mbox{{\rm \O}}}

\def\qed{\hfill \rule{4pt}{7pt}}

\def\pf{\noindent {\it Proof. }}

\begin{document}
	
	\title{Structure and coloring of some ($P_7$, $C_4$)-free graphs}
	\author{Ran Chen$^{1,}$\footnote{Email: 1918549795@qq.com},  \; Di Wu$^{1,}$\footnote{Email: 1975335772@qq.com}, \; Baogang  Xu$^{1,}$\footnote{Email: baogxu@njnu.edu.cn OR baogxu@hotmail.com. Supported by NSFC 11931006}\\\\
		\small $^1$Institute of Mathematics, School of Mathematical Sciences\\
		\small Nanjing Normal University, 1 Wenyuan Road,  Nanjing, 210023,  China}
	\date{}
	
	\maketitle
\begin{abstract}
Let $G$ be a graph. We use $P_t$ and $C_t$ to denote a path and a cycle on $t$ vertices, respectively. A {\em diamond} is a graph obtained from two triangles that share exactly one edge. A {\em kite} is a graph consists of a diamond and another vertex adjacent to a vertex of degree 2 of the diamond. A {\em gem} is a graph that consists of a $P_4$ plus a vertex adjacent to all vertices of the $P_4$. In this paper, we prove some structural properties to $(P_7, C_4,$ diamond)-free graphs, $(P_7, C_4,$ kite)-free graphs and $(P_7, C_4,$ gem)-free graphs. As their corollaries, we show that (\romannumeral 1) $\chi (G)\leq \max\{3,\omega(G)\}$ if $G$ is $(P_7, C_4,$ diamond)-free, (\romannumeral 2) $\chi(G)\leq \omega(G)+1$ if $G$ is $(P_7, C_4,$ kite)-free and (\romannumeral 3) $\chi(G)\leq 2\omega(G)-1$ if $G$ is $(P_7, C_4,$ gem)-free. These conclusions generalize some results of Choudum {\em et al} \cite{CKB21} and Lan {\em et al} \cite{LZL22}.

\begin{flushleft}
	{\em Key words and phrases:} $P_7$-free graphs, chromatic number, clique number\\
	{\em AMS 2000 Subject Classifications:}  05C15, 05C75\\
\end{flushleft}

\end{abstract}

\newpage

\section{Introduction}
All graphs considered in this paper are finite and simple. We follow \cite{BM08} for undefined notations and terminology. We use $P_t$ and $C_t$ to denote a path and a cycle on $t$ vertices, respectively. Let $ G $ be a graph, let $v\in V(G)$, and let $X$ and $Y$ be two subsets of $V(G)$. We say that $v$ is {\em complete} to $X$ if $v$ is adjacent to all vertices of $X$, and say that $v$ is {\em anticomplete} to $X$ if $v$ is not adjacent to any vertex of $X$. We say that $X$ is complete (resp. anticomplete) to $Y$ if each vertex of $X$ is complete (resp. anticomplete) to $Y$.

A {\em hole} of $G$ is an induced cycle of length at least 4, and a {\em $k$-hole} is a hole of length $k$.
Given a graph $H$, we say that $G$ is $H$-free if it dose not contain an induced subgraph isomorphic to $H$. Given a set $\{H_1, H_2, \cdots\}$ of graphs, we say that $G$ is $(H_1, H_2, \cdots)$-free if $G$ is $H_i$-free for each $i$.

A {\em clique blowup} of a graph $H$ with  $V(H)=\{v_1, v_2, \cdots, v_n\}$ is any graph $G$ such that $V(G)$ can be partitioned into $n$ nonempty cliques $A_i$, $v_i\in V(H)$ such that $A_i$ is complete to $A_j$ if $v_i\sim v_j $, and $A_i$ is anticomplete to $A_j$ if $v_i\not\sim v_j$. Particularly, if $|A_i|=t$ for all $i\in\{1, \cdots, n\}$, then we call $G$ a $t$-{\em size clique blowup} of $H$.

For $v\in V(G)$, let $N_G(v)$ be the set of vertices adjacent to $v$, $d_G(v)=|N_G(v)|$, $N_G[v]=N_G(v)\cup \{v\}$, and $M_G(v)=V(G)\setminus N_G[v]$. For $ X\subseteq V(G) $, let $N_G(X)=\{u\in V(G)\setminus X\;|\; u$ has a neighbor in $X\}$ and $M_G(X)=V(G)\setminus (X\cup N_G(X))$. Moreover,  $ G[X] $ denotes the subgraph of $ G $ induced by $ X $. If it does not cause any confusion, we usually omit the subscript $G$ and simply write $N(v)$, $d(v)$, $ N[v]$, $M(v)$, $N(X)$ and $M(X)$.  Let $\delta(G)$ denote the minimum degree of $G$. For $u$, $v\in V(G)$, we simply write $u\sim v$ if $uv\in E(G)$, and write $u\not\sim v$ if $uv\not\in E(G)$.

A {\em clique} ({\em stable~set}) of $ G $ is a set of mutually adjacent (non-adjacent) vertices in $ G. $ The clique number of $G$, denoted by $ \omega(G), $ is the maximum size of a clique in $G$. For positive integer $k$, a $k$-{\em coloring } of $G$ is a function $c: V(G)\rightarrow \{1,\cdots,k\}$, such that for each edge $uv$ of $G$, $c(u)\not=c(v)$. The chromatic number of $G$, denoted by $\chi(G)$, is the minimum number $k$ for which there exists a $k$-coloring of $G$. A graph is {\em perfect} if all its induced subgraphs $H$ satisfy $\chi(H)=\omega(H)$. A family ${\cal G}$ of graphs is said to be $\chi$-{\em bounded} \cite{G75} if there exists a function $f$ such that for every graph $G\in {\cal G}$, $\chi(G)\leq f(w(G))$. If such a function $f$ does exist to ${\cal G}$, then $f$ is called a $\chi$-{\em binding~function} of ${\cal G}$.

A clique cutset of $G$ is a clique $K$ in $G$ such that $G-K$ has more components than $G$. A vertex $u$ is a cut vertex of $G$ if it is a clique cutset of size 1.

\begin{figure}[htbp]\label{fig-1}
	\begin{center}
		\includegraphics[width=8cm]{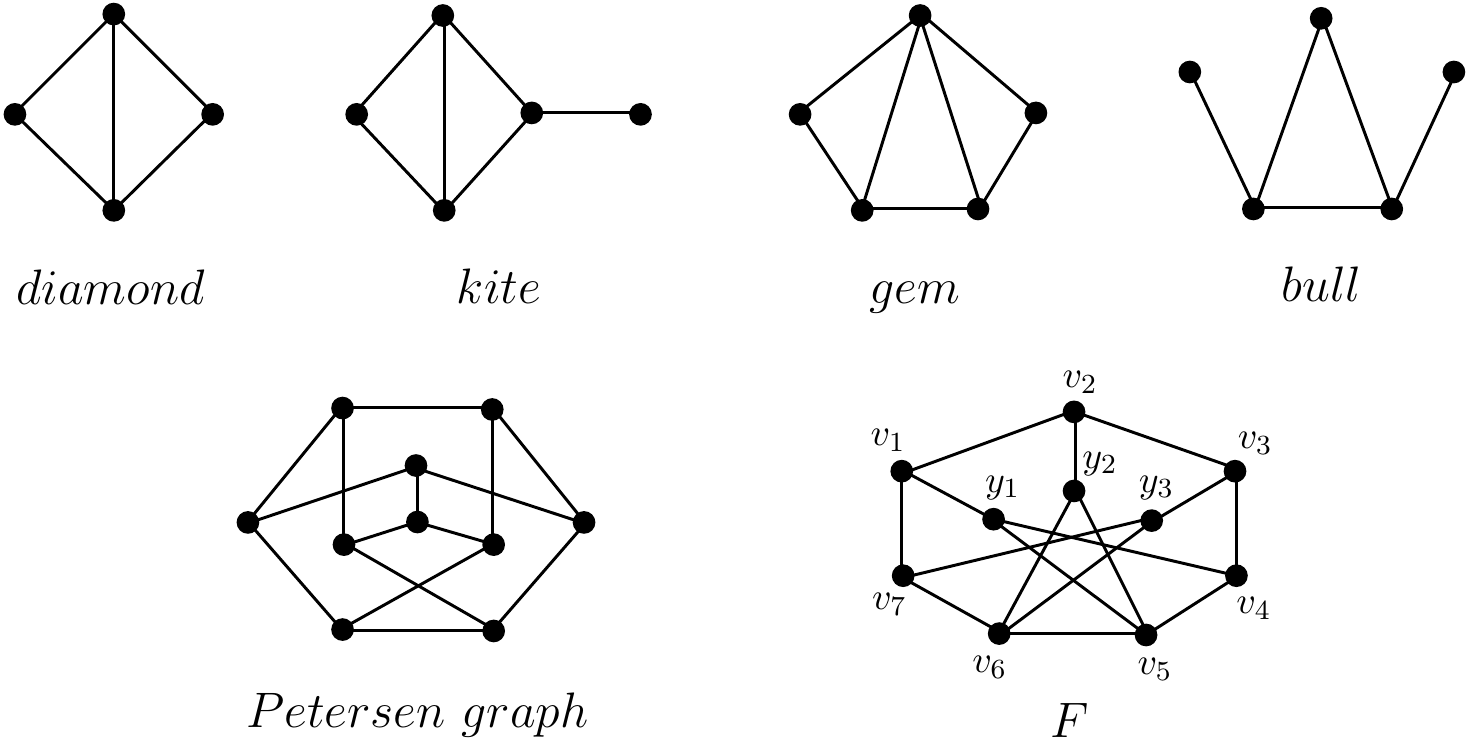}
	\end{center}
	\vskip -25pt
	\caption{Illustration of diamond, kite, gem, bull, Petersen graph and $F$.}
\end{figure}

A {\em diamond} is a graph obtained from two triangles that share exactly one edge. A {\em kite} is a graph consists of a diamond and another vertex adjacent to a vertex of degree 2 of the diamond. A {\em gem} is a graph that consists of a $P_4$ plus a vertex adjacent to all vertices of the $P_4$. A {\em bull} is a graph consisting of a triangle with two disjoint pendant edges. Let $C=v_1v_2v_3v_4v_5v_6v_7v_1$ be a 7-hole. We use $F$ to denote a graph obtained from $C$ by adding a stable set $\{y_1, y_2, y_3\}$ such that $N(y_1)\cap V(C)=\{v_1, v_4, v_5\}$, $N(y_2)\cap V(C)=\{v_2, v_5, v_6\}$ and $N(y_3)\cap V(C)=\{v_3, v_6, v_7\}$. (See Figure \ref{fig-1}).

Gy\'{a}rf\'{a}s \cite{G75} showed that $\chi(G)\leq (t-1)^{\omega(G)-1}$  for all $P_t$-free graphs. This upper bound was improved to $\chi(G)\leq (t-2)^{\omega(G)-1}$ in \cite{GHM03}. However, this $\chi$-binding function is exponential in $\omega(G)$. It is natural to ask that is it possible to improve the exponential bound for $P_t$-free graphs to a polynomial bound \cite{S16}. Interested readers are referred to \cite{RS2004, SR2019, SS18} for more results and still open problems on this topic.

It is known that $P_3$-free graphs are complete graphs and $P_4$-free graphs are  perfect  \cite{S74}. But up to now, no polynomial binding function for $P_t$-free graphs has been found, when $t\geq 5$. We refer the interested readers to \cite{G22}  for results of $P_5$-free graphs, and list here some conclusions on binding functions of some $P_6$-free or $ P_7$-free graphs.

Choudum, Karthick and Shalu \cite{CKS07} proved that every ($P_6$, gem)-free graph $G$ satisfies $\chi(G)\leq 8\omega(G)$. Cameron, Huang and Merkel \cite{CHM18} proved that every ($P_6$, diamond)-free graph $G$ satisfies $\chi(G)\leq \omega(G)+3$. Gasper and Huang \cite{GSH17} proved that every $(P_6, C_4)$-free graph $G$ satisfies $\chi(G)\leq \lfloor\frac{3}{2}\omega(G)\rfloor$. This bound was recently improved by Karthick and Maffray \cite{KM19} to $\chi (G)\leq\lceil\frac{5}{4}\omega(G) \rceil$. Mishra \cite{M21} proved that every ($P_6$, diamond, bull)-free graphs $G$ satisfies $\chi(G)\leq4$ when $\omega(G)=2$ and $\chi(G)=\omega(G)$ when $\omega(G)\ne 2$, and every ($P_7$, diamond, bull)-free graph $G$ satisfies $\chi(G)\leq\max\{7, \omega(G)\}$. Cameron, Huang, Penev and Sivaraman \cite{CHPS20} showed that every $(P_7, C_4, C_5)$-free graph satisfies $\chi(G)\leq\frac{3}{2}\omega(G)$. This upper bound was recently improved by Huang \cite{H22} to $\chi(G)\leq\lceil\frac{11}{9}\omega(G)\rceil$.

In 2021, Choudum, Karthick and Belavadi \cite{CKB21} proved that $\chi(G)\le 2\omega(G)-1$ if $G$ is $(P_7, C_7, C_4$, gem)-free, and
$\chi(G)\le \max\{3, \omega(G)\}$ if $G$ is $(P_7, C_7, C_4$, diamond)-free.  In 2022, Lan, Zhou and Liu \cite{LZL22} proved that $\chi(G)\le \max\{3, \omega(G)\}$ if $G$ is $(P_6, C_4$, diamond)-free.

In this paper, we study ($P_7, C_4, H)$-free graphs for $H\in$ \{diamond, kite, gem\}, and get some structural properties of these graphs. For graphs $G$ and $H$, we use $G+H$ to denote a graph with vertex set $V(G)\cup V(H)$ and edge set $E(G)\cup E(H)\cup \{uv\;|\; u\in V(G), v\in V(H)\}$. A {\em bisimplicial} vertex is one whose neighborhoods is the union of two cliques.

\begin{theorem}\label{diamond}
	Let $G$ be a connected $(P_7, C_4$, diamond)-free graph without clique cutsets. If $G$ is isomorphic to neither $F$ nor the Petersen graph, then $\delta(G)\le\max \{2, \omega(G)-1\}$.
\end{theorem}

\begin{theorem}\label{kite}
	Let $G$ be a connected $(P_7, C_4$, kite)-free graph with $\delta(G)\ge \omega(G)+1$. If $G$ has no clique cutsets, then there is an integer $\ell\ge 0$ such that $G=K_{\ell}+H$, where $H$ is the Petersen graph or $F$.	
\end{theorem}

\begin{theorem}\label{gem}
	Let $G$ be a connected $(P_7, C_4$, gem)-free graph without clique cutsets. If $G$ is not the clique blowup of the Petersen graph, then $G$ has a bisimplicial vertex.
\end{theorem}

There exist graphs showing that all the requirements in our theorems are necessary. We will give these examples separately after completing their proofs.  As immediate consequences of these theorems, we have the following corollaries that generalize some results of Choudum, Karthick and Belavadi \cite{CKB21}, and of Lan, Zhou and Liu \cite{LZL22}.

\begin{corollary}\label{diamond*}
	Every  $(P_7, C_4$, diamond)-free graph $G$ satisfies $\chi (G)\leq \max\{3,\omega(G)\}$.
\end{corollary}

\begin{corollary}\label{kite*}
	Every  $(P_7, C_4$, kite)-free graph $G$ satisfies $\chi(G)\leq\omega(G)+1$.
\end{corollary}

\begin{corollary}\label{gem*}
	Every  $(P_7, C_4$, gem)-free graph $G$ satisfies $\chi (G)\leq 2\omega(G)-1$.
\end{corollary}

We will prove Theorem \ref{diamond} in Section 2, prove Theorem \ref{kite} in Section 3, and prove Theorem \ref{gem} in Section 4.

\section{The class of $(P_7, C_4,$ diamond)-free graphs}

In this section, we focus on $(P_7, C_4$, diamond)-free graphs, and prove Theorem~\ref{diamond} and Corollary~\ref{diamond*}. Theorem~\ref{diamond} generalizes the following result from \cite{CKB21}.

\begin{lemma}\label{diamond'}{\em\cite{CKB21}}
Let $G$ be a connected $(P_7, C_7, C_4$, diamond)-free graph. Then $G$ has a clique cutset or $\delta(G)\leq max \{2, \omega(G)-1\}$ or $G$ is the Petersen graph.
\end{lemma}

\medskip

\noindent{\em Proof of Theorem~\ref{diamond}}: Let $G$ be a connected $(P_7, C_4$, diamond)-free graph. If $G$ is $C_7$-free, then the statement follows directly from Lemma~\ref{diamond'}. So, we suppose that $G$ has 7-holes, and let $v_1v_2v_3v_4v_5v_6v_7v_1$ be a 7-hole of $G$. Let $A=\{v_1, v_2, \cdots, v_7\}$, and let $R=M(A)$. During the proof of Theorem~\ref{diamond}, every subscript is understood to be modulo 7. For each $i\in$\{1, $\cdots$, 7\}, let
\begin{eqnarray*}	
	X_i&=&\{x\in N(A)|~N(x)\cap A=\{v_i, v_{i+3}\}\},~\mbox{and}\\
	Y_i&=&\{x\in N(A)|~N(x)\cap A=\{v_i, v_{i+3}, v_{i+4}\}\}.
\end{eqnarray*}
Let $X=X_1\cup\cdots\cup X_7$ and $Y=Y_1\cup\cdots\cup Y_7$.  Next, we prove some useful properties of $N(A)$.

\medskip

\begin{equation}\label{NA-1}
	\mbox{$N(A)=X\cup Y$, and $V(G)=A\cup X\cup Y\cup R$.}
\end{equation}

 It is certain that $X\cup Y\subseteq N(A)$. Since $G$ is diamond-free, we have that no vertex of $N(A)$ may have three consecutive neighbors in $A$. Let $x\in N(A)$. By symmetry, we may assume $x\sim v_1$. If $N(x)\cap A=\{v_1\}$, then $G[\{x, v_1, v_2, v_3, v_4, v_5, v_6\}]=P_7$, a contradiction. Therefore $|N(x)\cap A|\geq2$.

Suppose $x\sim v_2$. Then, $x$ has neighbors in $\{v_3, \cdots, v_7\}$ as otherwise  $G[\{x, v_2, v_3, v_4, v_5, v_6, v_7\}]=P_7 $, and $x\not\sim v_3$ and $x\not\sim v_7$ as otherwise $x$ has three consecutive neighbors in $A$. To forbid a 4-hole on $\{x,v_2,v_3,v_4\}$ or $\{x,v_1,v_6,v_7\}$, we have that $x\not\sim x_4$ and $x\not\sim x_6$. Therefore $x\sim v_5$, which implies that $x\in Y$. By symmetry, we have $x\in Y$ when $x\sim v_7$. So, we way assume that $x\not\sim v_2$ and $x\not\sim v_7$, and thus $x\not\sim v_3$ and $x\not\sim v_6$ to forbid a 4-hole on $\{x, v_1, v_2, v_3\}$ or $\{x, v_1, v_6, v_7\}$.

Now, $x\in X$ if $x$ has exactly one neighbor in $\{v_4, v_5\}$, and $x\in Y$ if $x$ is complete to $\{v_4, v_5\}$. This proves (\ref{NA-1}).

\begin{equation}\label{NA-2}
\mbox{For each $i\in\{1,\cdots,7\}$, $|X_i|\leq 1$ and $ |Y_i|\leq 1$.}
\end{equation}

Suppose to its contrary, we may assume by the symmetry that $X_1$ or $ Y_1 $ has two distinct vertices $x$ and $x'$. Then $G[\{x, x', v_1, v_4\}]$ is a diamond if $x\sim x'$, and $xv_1x'v_4x$ is a 4-hole if $x\not\sim x'$. Both are contradictions. Therefore, (\ref{NA-2}) holds.

For $i\in\{1,\cdots,7\}$, let $X_i=\{x_i\}$ if $X_i\ne\emptyset$, and $Y_i=\{y_i\}$ if $Y_i\ne\emptyset$.

\begin{equation}\label{NA-3}
\mbox{$N(A)$ is a stable set.}
\end{equation}

Suppose to its contrary that $N(A)$ has two adjacent vertices $u$ and $v$. First we assume that $\{u, v\}\cap X\ne\emptyset.$ Without loss of generality, suppose that $ u\in X_1. $ Then $v\not\sim v_j$ for $j\in$\{2, 7\} as otherwise $G[\{u, v, v_1, v_j\}]$ is a diamond if $v\sim v_1$, and $uvv_jv_1u$ is a 4-hole if $v\not\sim v_1$. Moreover, $v\not\sim v_k$ for $k\in\{3, 5\}$ as otherwise $G[\{u, v, v_4, v_k\}]$ is a diamond if $v\sim v_4$, and $uvv_kv_4u$ is a 4-hole if $v\not\sim v_4$. Therefore $N(v)\cap A\subseteq\{v_1, v_4, v_6\}$. Since $|N(v)\cap A|\ge2$ and $|X_1|\leq 1$, we have that $N(v)\cap A\in\{\{v_1, v_6\}, \{v_4, v_6\}, \{v_1, v_4, v_6\}\}$. Then $vv_6v_7v_1v$ is a 4-hole if $N(v)\cap A\in\{\{v_1, v_6\}, \{v_1, v_4, v_6\}\}$, and $vv_4v_5v_6v$ is a 4-hole if $N(v)\cap A=\{v_4, v_6\}$, both are contradictions.

So, we may assume that $u\in Y$ and $v\in Y$. Without loss of generality, suppose that $u\in Y_1$. Then $v\not\sim v_j$ for $j\in\{2, 7\}$ as otherwise $G[\{u,v,v_1,v_j\}]$ is a diamond if $v\sim v_1$, and $uvv_jv_1u$ is a 4-hole if $ v\not\sim v_1$. Moreover, $v\not\sim v_3$ to forbid a diamond on $\{u,v,v_3,v_4\}$ if $v\sim v_4$ or a 4-hole on $\{u,v,v_3,v_4\}$ if $v\not\sim v_4$. If $v\sim v_5$, then $v\not\sim v_4$ as otherwise $v\in Y_1$ and thus $|Y_1|\geq 2$, a contradiction. But now $G[\{u,v,v_4,v_5\}]$ is a diamond. Therefore, $N(v)\cap A=\{v_1,v_4,v_6\}$ since $|N(v)\cap A|=3$. However, $vv_4v_5v_6v$ is a 4-hole, a contradiction. This proves (\ref{NA-3}).

\begin{equation}\label{NA-4}
	\mbox{For each $i\in\{1,\cdots,7\}$, ethier $X_i=\emptyset$ or $X_{i+2}\cup X_{i+5}=\emptyset$.}
\end{equation}

 If this is not true, we may assume by symmetry that $X_1\ne\emptyset$ and $X_3\ne\emptyset$, then $x_1\not\sim x_3$ by (\ref{NA-3}), which implies that $G[\{x_1, x_3, v_1, v_2, v_4, v_5, v_6\}]=P_7$, a contradiction. This proves (\ref{NA-4}).

\begin{equation}\label{NA-5}
	\mbox{For each $i\in\{1, \cdots, 7\}$, ethier $Y_i=\emptyset$ or $Y_{i+3}\cup Y_{i+4}=\emptyset$.}
\end{equation}

 If (\ref{NA-5}) does not hold, we may assume by symmetry that $Y_1\ne\emptyset$ and $Y_4\ne\emptyset$, then $y_1\not\sim y_4$  by (\ref{NA-3}), which implies that $y_1v_4y_4v_1y_1$ is a 4-hole, a contradiction. This proves (\ref{NA-5}).

\begin{equation}\label{NA-6}
\mbox{ For each $i\in\{1, \cdots, 7\}$, ethier $X_i=\emptyset$ or $Y_{i}\cup Y_{i+1}\cup Y_{i+2}\cup Y_{i+3}=\emptyset$.}
\end{equation}

 If it is not the case, we may assume by symmetry that $X_1\ne\emptyset$ and $Y_1\cup Y_2\cup Y_3\cup Y_4\ne\emptyset$.  Let $y\in Y_1\cup Y_2\cup Y_3\cup Y_4$. Then $x_1\not\sim y$ by (\ref{NA-3}), which implies that $x_1v_4yv_1x_1$ is a 4-hole if $y\in Y_1\cup Y_4$, $G[\{x_1, y, v_2, v_3, v_4, v_6, v_7\}]=P_7$ if $y\in Y_2$, and $G[\{x_1, y, v_1, v_2, v_4, v_5, v_6\}]=P_7$ if $y\in Y_3$, all of which are contradictions. This proves (\ref{NA-6}).

\medskip

Suppose that $Y=\emptyset.$ By (\ref{NA-4}), we have that $\{X_1, X_2, \cdots, X_7\}$ has at most three nonempty elements, and thus $d(v_i)=2$ for some $i\in\{1, \cdots, 7\}$. If $\omega(G)=2$, then $Y=\emptyset$, and thus $\delta(G)\leq 2$.

Suppose that $\omega(G)\geq 4$ and suppose by symmetry that $Y_1\ne \emptyset$. Then, $Y_4\cup Y_5\cup X_1\cup X_7\cup X_6\cup X_5=\emptyset$ by (\ref{NA-5}) and (\ref{NA-6}), and one can  verify easily that $d(v_1)=3\leq \omega(G)-1\leq \max\{2, \omega(G)-1\}$.

So, we suppose that
\begin{itemize}
\item $\omega(G)=3$ and $\delta(G)\geq 3$,
\item $G$ is not isomorphic to $F$ and has no clique cutsets, and
\item $Y_{i_0}=\{y_{i_0}\}\neq\emptyset$ for some integer $i_0$.
\end{itemize}

By (\ref{NA-5}) and (\ref{NA-6}), we have that  $Y_{i_0+3}\cup Y_{i_0+4}\cup X_{i_0}\cup X_{i_0-1}\cup X_{i_0-2}\cup X_{i_0-3}=\emptyset$, and so
\begin{equation}\label{eqa-X-Y-0}
X=X_{i_0+1}\cup X_{i_0+2}\cup X_{i_0+3}\mbox{ and } Y=Y_{i_0}\cup Y_{i_0+1}\cup Y_{i_0+2}\cup Y_{i_0-2}\cup Y_{i_0-1}.
\end{equation}
We will show that
\begin{equation}\label{eqa-ca1-1}
	\mbox{$X_{i_0+1}=\emptyset$ if $Y_{i_0}\neq\emptyset$}.
\end{equation}

For simplicity, we may take $i_0=1$ by symmetry to illustrate the procedure. Before proving (\ref{eqa-ca1-1}), we firstly prove that

\begin{equation}\label{ca1-1}
	\mbox{$X_{i_0+2}=\emptyset$ if $X_{i_0+1}\ne\emptyset$}.
\end{equation}

Notice that we take $i_0=1$. Suppose $X_2\ne\emptyset$ and $X_3\ne\emptyset$. By (\ref{NA-4}) and (\ref{NA-6}), we have that $X_4\cup X_7\cup Y_2\cup Y_3\cup Y_4\cup Y_5\cup Y_6=\emptyset$, and so $X=X_2\cup X_3$ and $Y=Y_1\cup Y_7$. Since $Y_7=\emptyset$ implies that $d(v_7)=2$, we have that $Y_7=\{y_7\}\neq\emptyset$.

Since $\delta (G)\geq 3 $ and $N(A)$ is a stable set, both $x_2$ and $x_3$ have neighbors in $R$. Let $u$ be a neighbor of $x_2$ in $R$ and $v$ be a neighbor of $x_3$ in $R$, and let $S$ and $T$ be the components of $G[R]$ which contains $u$ and $v$, respectively. Since $G[\{x_2, v_1, v_2, v_6, v_7\}]=P_5$ and $G$ is $P_7$-free, we have that $x_2$ must be complete to $V(S)$, and thus $\omega(S)\le2 $ as $\omega(G)=3$. Moreover, $S$ is $P_3$-free as $G$ is diamond-free, and thus $S$ is a $K_1$ or a $K_2$. Similarly, $x_3$ must be complete to $V(T)$, and $T$ is a $K_1$ or a $K_2$.

Suppose that $S=T$. Then $|V(S)|=|V(T)|=1$ as otherwise $G[\{x, y, x_2, x_3\}]$ is a diamond for any edge $xy$ of $S$, and so $V(S)=V(T)=\{u\}=\{v\}$. If $u\sim y_7$, then $v_3x_3uy_7v_3$ is a 4-hole, a contradiction. So, $u\not\sim y_7$. By symmetry, $u\not\sim y_1$. But now $d(u)=2$, a contradiction. Therefore, $S\ne T$.

If $V(S)=\{u\}$, then we have that $u\not\sim y_1$ to forbid a 4-hole on $\{u, y_1, x_2, v_5\}$, and thus $u\sim x_3$ and $u\sim y_7$ since $\delta (G)\geq 3$, which implies that there is a 4-hole on $\{u, x_3, v_3, y_7\}$, a contradiction. So, $S$ is a $K_2$. Similarly, $T$ is a $K_2$. Let $V(S)=\{u, u'\}$ and $V(T)=\{v, v'\}$.

To forbid a 4-hole on $\{y_1, x_2, v_5, u\}$ or $\{y_1, x_2, v_5, u'\}$, $y_1$ must be anticomplete to $V(S)$. Similarly, $y_7$ must be anticomplete to $V(T)$. Since $\delta(G)\geq 3$, we have that each vertex of $S$ has neighbors in $\{x_3, y_7\}$. Since $G$ is diamond-free, we may assume, without loss of generality, that $u\sim x_3$ and $u'\sim y_7 $. Similarly, we may suppose that $v\sim x_2$ and $v'\sim y_1$. But now, $ux_2vx_3u$ is a 4-hole, a contradiction. This proves (\ref{ca1-1}).

\medskip

We can now turn to prove (\ref{eqa-ca1-1}). Recall that we take $i_0=1$ by symmetry.

Suppose to its contrary that $X_2\ne\emptyset $. Then $X_3=\emptyset$ by (\ref{ca1-1}). By (\ref{NA-4}) and (\ref{NA-6}), we have that $ X_4\cup X_7\cup Y_2\cup Y_3\cup Y_4\cup Y_5=\emptyset$, and thus $X=X_2$ and $Y=Y_1\cup Y_6\cup Y_7$. Since, for $i\in\{6, 7\}$, $Y_i=\emptyset$ implies that $d(v_i)=2$, we have that $Y_6\neq\emptyset$ and $Y_7\neq\emptyset$. Since $\delta(G)\geq 3$ and $N(A)$ is a stable set,  we have that $x_2$ must have a neighbor, say $u$, in $R$. Let $B$ be the component of $G[R]$ that contains $u$.

Since $G$ is $P_7$-free and diamond-free, we have that $x_2$ must be complete to $V(B)$ and $B$ is a $K_1$ or a $K_2$. Suppose that $V(B)=\{u\}$. Then, $u$ has at least two neighbors in $N(A)\setminus \{x_2\}$ as $\delta(G)\ge3$. Since $u\sim y_1$ implies that there is a 4-hole on $\{u, y_1, x_2, v_5\}$, we have that $u\not\sim y_1$, and thus $u\sim y_6$ and $u\sim y_7 $. But then, $uy_6v_3y_7u$ is a 4-hole, a contradiction. So $B$ is a $K_2$. Let $V(B)=\{u, v\}$. To forbid a 4-hole on $\{a, y_1, x_2, v_5\}$ or $\{a, x_2, v_2, y_6\}$, where $a\in V(B)$, we have that $V(B)$ is anticomplete to \{$y_1$, $y_6$\}. Since $\delta(G)\geq 3$, each vertex in $B$ must have neighbor in $\{y_7\}$, which implies a diamond $G[\{u, v, x_2, y_7\}]$, a contradiction. This proves (\ref{eqa-ca1-1}).

\medskip

Next, we prove that
\begin{equation}\label{eqa-emptyset-X}
X=\emptyset \mbox{ and } Y=Y_{i_0}\cup Y_{i_0+1}\cup Y_{i_0+2}\cup Y_{i_0-1}\cup Y_{i_0-2}.
\end{equation}

We still take $i_0=1$ for simplicity. By (\ref{eqa-X-Y-0}) and (\ref{eqa-ca1-1}), we have that $X=X_3\cup X_4$ and $Y=Y_1\cup Y_2\cup Y_3\cup Y_6\cup Y_7$.

Suppose $X_3\neq\emptyset$. By (\ref{NA-6}), we have that $Y_3\cup Y_4\cup Y_5\cup Y_6=\emptyset$, and so $X=X_3\cup X_4$ and $Y=Y_1\cup Y_2\cup Y_7$. Since $d(v_2)\ge 3$, we have that $Y_2\neq\emptyset$, and consequently $X_3=\emptyset$ by taking $i_0=2$ in (\ref{eqa-ca1-1}), a contradiction. Therefore,  $X_3=\emptyset$, and $X=X_4$ and $Y=Y_1\cup Y_2\cup Y_3\cup Y_6\cup Y_7$.

Suppose that $X_4\neq\emptyset$. Then, $Y_4\cup Y_5\cup Y_6\cup Y_7=\emptyset$ by (\ref{NA-6}), and so $Y=Y_1\cup Y_2\cup Y_3$. Moreover, $Y_2\neq\emptyset$ and $Y_3\neq\emptyset$ as $d(v_2)\ge 3$ and $d(v_3)\ge 3$. Consequently, we have a contradiction as $X_4=\emptyset$ by taking $i_0=3$ in  (\ref{eqa-ca1-1}). Therefore, $X_4=\emptyset$, and thus (\ref{eqa-emptyset-X}) holds.

\medskip

Recall that $Y_{i_0}\ne \emptyset$ by our assumption. By symmetry, we may set $i_0=1$ in (\ref{eqa-emptyset-X}). Then, we have that
\begin{equation}\label{eqa-X-Y12367}
\mbox{$X=\emptyset$ and $Y=Y_{1}\cup Y_{2}\cup Y_{3}\cup Y_{6}\cup Y_{7}$.}
\end{equation}

We prove further that
\begin{equation}\label{ca1-3}
	\mbox{$R=\emptyset$ if $Y_2\ne\emptyset$ and $Y_3\ne\emptyset$}.
\end{equation}

Suppose that $Y_2\ne\emptyset$, $Y_3\ne\emptyset$, and $R\ne\emptyset$. By (\ref{NA-5}), we have that $Y_5\cup Y_6\cup Y_7=\emptyset$, and thus $Y=Y_1\cup Y_2\cup Y_3=\{y_1, y_2, y_3\}$. Let $Q$ be a component of $G[R]$. Since $G$ is connected, there must exist a path $P$ between $Q$ and $Y$. If we can show that $P$ must contain $y_2$, then $y_2$ is a cutvertex of $G$, which leads to a contradiction. Suppose that $V(P)\cap Y=\{y\}$. To prove (\ref{ca1-3}), it suffices to prove the following Claim~\ref{c-1}.

\begin{claim}\label{c-1}
$y\ne y_3$ and $y\ne y_1. $
\end{claim}
\pf By symmetry, we only need to prove $ y\ne y_3 $. Suppose to its contrary that $y=y_3 $. Since $G[\{v_1, v_2, v_3, y_3\}]=P_4$ and $G$ is $P_7$-free, we have that $d(u, y_3)\le2$ for each vertex $u$ of $Q$. Let $N_i=\{u\in V(Q)|~d(u,y_3)=i\}$ for $i\in$\{1, 2\}.  Since $G$ is diamond-free and $y_3$ is complete to $N_1$, and since $\omega(G)=3$, we have that each component of $G[N_1]$ is a $K_1$ or a $K_2$. We show that
\begin{equation}\label{ca1-4}
	\mbox{$N_2$ is a stable set}.
\end{equation}

Suppose to the contrary that $N_2$ has two adjacent vertices $x'$ and $y'$. Let $S'$ be the component of $G[N_2]$ which contains $x'$ and $y'$. Let $ x''$ be a neighbor of $x'$ in $ N_1$. Since $d(u, y_3)\le2$ for each vertex $u$ of $Q$, we have that $x''$ is complete to $V(S')$, which implies that $\omega(S')=2$ as $\omega(G)=3$. Moreover, $S'$ is $P_3$-free as $G$ is diamond-free. Then, $S'=x'y'$.

Let $S''$ be the component of $G[N_1]$ which contains $x''$. Note that $S''$ is a $K_1$ or a $K_2$. For each vertex $t_1\in N_1\setminus\{x''\}$ and $t_2\in V(S')$, to avoid a 4-hole on $\{x'',y_3,t_1,t_2\}$ if $t_1\not\sim x''$ or a diamond on $\{y_3,x'',t_1,t_2\}$ if $t_1\sim x''$, $\{x', y'\}$ must be anticomplete to $N_1\setminus\{x''\}$, which implies that $x'$ and $y'$ both have a neighbor in $\{y_1, y_2\}$ as $\delta(G)\geq 3$.

Suppose $x'\sim y_1$ and $y'\sim y_1$. Then $x''\sim y_1$ to forbid a diamond on $\{x', y', x'', y_1\}$. If $y_2$ has a neighbor in $S'$, say $x'$ by symmetry, then $y_2\sim y'$ to forbid an induced $P_7$ on $\{x', y', y_2, y_3, v_1, v_2, v_7\}$, and thus $G[\{x', y', y_1, y_2\}]$ is a diamond, a contradiction. So, $y_2$ has no neighbors in $S'$, and thus $\{x'', y_1\}$ is a clique cutset, which leads to  a contradiction. This proves that $x'\not\sim y_1$ or $y'\not\sim y_1$.

If $x'\sim y_2$ and $y'\sim y_2$, then $x''\not\sim y_2$ to forbid a 4-hole on $\{x'', y_2, y_3, v_6\}$, which forces $G[\{x', y', x'', y_2\}]$ to be a diamond, a contradiction. Therefore, $x'\not\sim y_2$ or $y'\not\sim y_2$.

By symmetry, we may assume that $x'\sim y_1$, $x'\not\sim y_2$, $y'\sim y_2$, and $y'\not\sim y_1$. Then, $G$ has an induced $P_7$ on $\{x', y', y_1, y_3, v_1, v_2, v_3\}$, a contradiction. This proves (\ref{ca1-4}).

Suppose that $N_2\not=\emptyset$. Let $w\in N_2$ and $w'$ be a neighbor of $w$ in $N_1$. To avoid a diamond on $\{w,w',w'',y_3\}$ if $w''\sim w'$ or a 4-hole on $\{w,w',w'',y_3\}$ if $w''\not\sim w'$ for each vertex $w''\in N_1\setminus\{w'\}$, $w$ must be anticomplete to $V(Q)\setminus\{w', w\}$, and so $w$ is complete to $\{y_1, y_2\}$ because $\delta(G)\geq 3$. But then, $wy_1v_5y_2w$ is a 4-hole, which is a contradiction. This shows that $N_2=\emptyset$, and thus $Q=N_1$.

Suppose that $N_1\neq\emptyset$. Let $M$ be a component of $G[N_1]$, and let $a\in V(M)$. To avoid a 4-hole $ay_3y_6y_2a$, we have that $a\not\sim y_2$. Since $\delta(G)\ge 3$, we have that $M$ is a $K_2$ of which both vertices are adjacent to $y_1$, which forces a diamond $G[V(M)\cup \{y_1, y_3\}]$ and leads to a contradiction. Therefore, $N_1=\emptyset$, and so $V(Q)=\emptyset$. This completes the proof of Claim~\ref{c-1}, and proves (\ref{ca1-3}). \qed

\medskip

If $Y_2\ne\emptyset$ and $Y_3\ne\emptyset$, then $G$ is isomorphic to $F$  by (\ref{eqa-X-Y12367}) and (\ref{ca1-3}).

If $Y_2=\emptyset$ and $Y_3\not=\emptyset$, then $Y_6\cup Y_7=\emptyset$ by (\ref{NA-5}), which implies that $Y=Y_1\cup Y_3$ and $d(v_2)=2$, a contradiction.

So, we suppose that $Y_3=\emptyset$ and $Y_2\ne\emptyset$. Then,  $Y_5=Y_6=\emptyset$ by (\ref{NA-5}). Since $Y_7=\emptyset$ implies that $d(v_7)=2$, we have that $Y_7\not=\emptyset$, and so $Y=Y_1\cup Y_2\cup Y_7=\{y_1, y_2, y_7\}$. By setting $i_0=7$ and replacing the tuple $(Y_1, Y_2, Y_3)$ with $(Y_7, Y_1, Y_2)$, we can deduce, with the same argument as that used in proving  (\ref{ca1-3}), that $R=\emptyset$ and thus $G$ is isomorphic to $F$. This totally completes the proof of Theorem~\ref{diamond}. \qed

\medskip

The following three classes of graphs show that the requirements $P_7$-free, $C_4$-free, and diamond-free are all necessary.

Let $G_1$ be a $t$-size clique blowup of a 7-hole with $ t\geq 2 $. It is certain that $G_1$ is $(P_7, C_4$)-free. Moreover, we can easily verify that $G_1$ does not satisfy the conclusion of Theorem~\ref{diamond}.

Let $G_2$ be a graph whose vertex-set is partitioned into seven stable sets $S_1$, $S_2$, $\cdots$, $S_7$ such that for each $i\in \{1, 2, \ldots, 7\}$, $|S_i|\geq 2$, and $S_i$ is complete to $S_{i+1}\cup S_{i-1}$ and anticomplete to $S_{i+2}\cup S_{i-2}\cup S_{i+3}\cup S_{i-3}$. It is certain that $G_2$ is ($P_7$, diamond)-free and does not satisfy the conclusion of Theorem~\ref{diamond}.

\begin{figure}[htbp]\label{fig-2}
	\begin{center}
		\includegraphics[width=3.5cm]{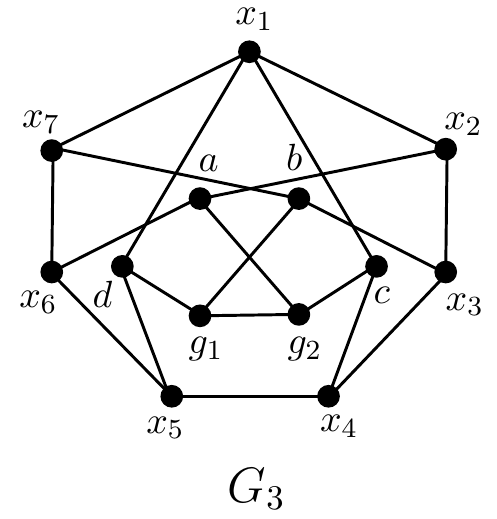}
	\end{center}
	\vskip -25pt
	\caption{Illustration of $G_3$.}
\end{figure}

Let $A_1=x_1x_2x_3x_4x_5x_6x_7x_1 $ be a 7-hole. We use $G_3$ to denote the graph obtained from $A_1$ by adding a tree with vertex set $\{a, b, c, d, g_1, g_2\}$ and edge set $\{ag_2, bg_1, cg_2, dg_1, g_1g_2\}$ such that $ N(a)\cap V(A_1)=\{x_2, x_6\}$, $N(b)\cap V(A_1)=\{x_3, x_7\}$, $N(c)\cap V(A_1)=\{x_1, x_4\}$, $N(d)\cap V(A_1)=\{x_1, x_5\}$ and $N(\{g_1, g_2\})\cap V(A_1)=\emptyset$. See Figure \ref{fig-2}. It is certain that $G_3$ is ($C_4$, diamond)-free, and does not satisfy the conclusion of Thoerem~\ref{diamond}.

Now, we can prove Corollary~\ref{diamond*} by induction on $|V(G)|$. Let $G$ be a $(P_7, C_4$, diamond)-free graph. We may assume that $G$ is connected and has no clique cutsets. By Theorem~\ref{diamond}, $\delta(G)\leq \max\{2, \omega(G)-1\}$, or $G$ is the Petersen graph or the graph $F$. It follows immediately that $\chi(G)\le\max\{3, \omega(G)\}$.

\section{The class of $(P_7, C_4,$ kite)-free graphs.}

In this section, we consider $(P_7, C_4$, kite)-free graphs, and prove Theorem~\ref{kite} and Corollary~\ref{kite*}. We will need the following Lemma~\ref{kite'} from \cite{FHH18}.

\begin{lemma}\label{kite'}{\em\cite{FHH18}}
Let $G$ be a connected ($C_4$, kite)-free graph. If $G$ has no clique cutsets, then $G=K_{\ell}+G'$ for some diamond-free graph $G'$, where $\ell$ is a nonnegative integer.
\end{lemma}

\medskip

\noindent\emph{Proof of Theorem~\ref{kite}}: Let $G$ be a connected ($P_7$, $C_4$, kite)-free graph such that $G$ has no clique cutsets and $\delta(G)\ge \omega(G)+1$.

By Lemma~\ref{kite'}, there exist an integer $\ell\ge 0$ and a diamond-free graph $G'$ such that $G=K_{\ell}+G'$. We have that $G'$ has no clique cutsets; Otherwise if $G'$ has a clique cutset $Q$, then $ Q\cup V(K_{\ell}) $ is a clique cutset of $G$, a contradiction. Therefore, $G'$ is a $(P_7, C_4$, diamond)-free graph without clique cutsets. By Theorem~\ref{diamond}, we have that $G'$ is isomorphic to $F$ or the Petersen graph, or $\delta(G')\leq \max\{2, \omega(G')-1\}$.

If $G'$ is isomorphic to $F$ or the Petersen graph, then we are done. So, we suppose that $G'$ is isomorphic to neither $F$ nor the Petersen graph. Thus, $\delta(G')\leq \max\{2, \omega(G')-1\}$.
If $G'$ is a single vertex, then $G=K_{\ell+1}$ and $\delta(G)=\ell=\omega(G)-1$. If $\omega(G')=2$, then $\delta(G')\le 2$, and so $\delta(G)\le \delta(G')+\ell\le \ell+2=\omega(G)$. If
$\omega(G')\ge 3$, then $\delta(G')\leq \omega(G')-1$, and so $\delta(G)\le \delta(G')+\ell\le \omega(G')+\ell-1\le \omega(G)-1$. All contradict our assumption that $\delta(G)\ge \omega(G)+1$. This proves Theorem~\ref{kite}. \qed
\medskip

Notice that the graph $G_1$ defined in Section 2 is $(P_7, C_4$)-free, and the graph $G_2$ defined in Section 2 is ($P_7$, kite)-free. Moreover, one can easily check that both $G_1$ and $G_2$ do not satisfy the conclusion of Theorem~\ref{kite}.

Let $A_2$ be the graph obtained from a 8-cycle $y_1y_2y_3y_4y_5y_6y_7y_8y_1 $ by adding a 6-cycle $u_1u_2u_3u_4u_5u_6u_1$ such that the edge set between them is $\{y_1u_1, y_4u_2, y_5u_4, y_6u_6\}$. We use $G_4$ to denote the graph obtained from $A_2$ by adding  a stable set $\{t_1, t_2\}$ such that $N(t_1)\cap V(A_2)= \{y_2, y_7, u_5\}$ and $N(t_2)\cap V(A_2)=\{y_3, y_8, u_3\}$ (see Figure 3). Obviously, $G_4$ is ($C_4$, kite)-free and it does not satisfy the conclusion of Theorem~\ref{kite}.

These graphs show that the requirements $P_7$-free, $C_4$-free, and kite-free in Theorem~\ref{kite} are all necessary.

\begin{figure}[htbp]\label{fig-3}
 	\begin{center}
 		\includegraphics[width=4cm]{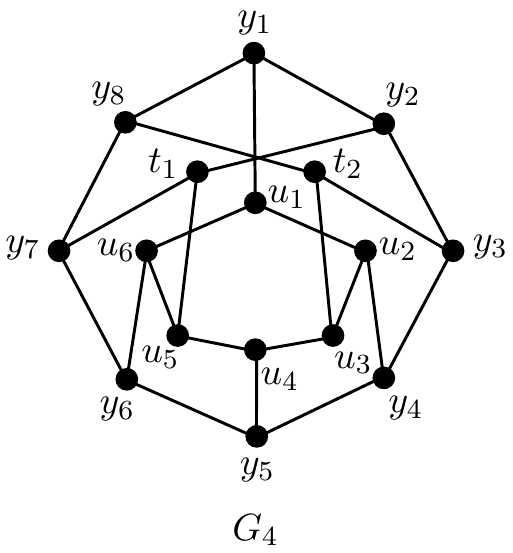}
 	\end{center}
 	\vskip -25pt
 	\caption{Illustration of $G_4$.}
\end{figure}

\medskip

Now, we can prove Corollary~\ref{kite*} by induction on $|V(G)|$. Let $G$ be a $(P_7, C_4$, kite)-free graph. We may assume that $G$ is connected and has no clique cutsets. By Theorem~\ref{kite}, we have that $\delta(G)\leq\omega(G)$ or there is an integer $\ell\geq0$ such that $G=K_\ell+H$, where $H$ is the Petersen graph or $F$. It follows immediately that $\chi(G)\le\omega(G)+1$.

\section{The class of $(P_7,C_4,$ gem)-free graphs.}

In this section, we turn our attention to $(P_7, C_4$, gem)-free graphs, and prove Theorem~\ref{gem} and Corollary~\ref{gem*}. Theorem~\ref{gem} generalizes the following Lemma~\ref{gem'} from \cite{CKB21}. In our proof, we will also use a conclusion from \cite{KM19}. Recall that a bisimplicial vertex is a vertex whose neighborhoods is the union of two cliques.

\begin{lemma}\label{gem'}{\em\cite{CKB21}}
	Let $G$ be a connected $(P_7, C_7, C_4$, gem)-free graph. If $G$ has no clique cutsets, then  $G$ has a bisimplicial vertex,  or $G$ is the clique blowup of the Petersen graph.
\end{lemma}

\begin{lemma}\label{Petersen}{\em\cite{KM19}}
If $G$ is any clique blowup of the Petersen graph, then $\chi(G)\leq \lceil \frac{5}{4}\omega(G) \rceil$.
\end{lemma}

\medskip

\noindent\emph{Proof of Theorem~\ref{gem}}: Let $G$ be a connected $(P_7, C_4$, gem)-free graph. If $G$ is $C_7$-free, then the statement follows directly from Lemma~\ref{gem'}. So, we suppose that $G$ has 7-holes, and let $v_1v_2v_3v_4v_5v_6v_7v_1$ be a 7-hole of $G$. Let $A=\{v_1,\cdots,v_7\}$ and let $R=M(A)$. During the proof of Theorem~\ref{gem}, every subscript is understood to be modulo 7.  Since $G$ is gem-free, it is certain that
\begin{equation}\label{eqa-no-four}
\mbox{no vertex of $N(A)$ may have four consecutive neighbors in $A$.}
\end{equation}

Let $i$ be an arbitrary integer in $\{1,\cdots,7\}$, let
\begin{eqnarray*}	
	X_i&=&\{x\in N(A)|~N(x)\cap A=\{v_i, v_{i+3}\}\},\\
	Y_i&=&\{x\in N(A)|~N(x)\cap A=\{v_i, v_{i+1}, v_{i+2}\}\},\\
	Z_i&=&\{x\in N(A)|~N(x)\cap A=\{v_i, v_{i+3}, v_{i+4}\}\},
\end{eqnarray*}
and let $X=\cup_{i=1}^7 X_i$, $Y=\cup_{i=1}^7 Y_i$ and  $Z=\cup_{i=1}^7  Z_i$. We will firstly prove some useful properties about $X_i$, $Y_i$ and $Z_i$.

\medskip

\noindent{\bf (M1)} $N(A)=X\cup Y\cup Z$, and $V(G)=A\cup X \cup Y\cup Z\cup R$.

It suffices to verify that $N(A)\subseteq X\cup Y\cup Z$.  Let $x\in N(A)$. Without loss of generality, we suppose that $x\sim v_1$. If $N(x)\cap A=\{v_1\}$, then $G[\{x, v_1, v_2, v_3, v_4, v_5, v_6\}]=P_7$, a contradiction. Therefore,  $|N(x)\cap A|\geq 2$.

Suppose $x\sim v_2$. To avoid an induced $P_7=xv_2v_3v_4v_5v_6v_7$, $x$ must have a neighbor in $\{v_3, \cdots, v_7\}$. If $x\sim v_3$, then $x\not\sim v_4$ and $x\not\sim v_7$ by (\ref{eqa-no-four}), and thus $x\not\sim v_5$ and $x\not\sim v_6$ to forbid a 4-hole $xv_3v_4v_5x$ or $xv_1v_7v_6x$, which leads to $x\in Y$ as $N(x)\cap A=\{v_1, v_2, v_3\}$. By symmetry, we have that $x\in Y$ if $x\sim v_7$. So, we may assume that $x\not\sim v_3$ and $x\not\sim v_7$. Consequently, $x\not\sim v_4$ and $x\not\sim v_6$ to forbid a 4-hole $xv_2v_3v_4x$ or $xv_1v_7v_6x$. Now, $x\sim v_5$ and hence $x\in Z$ as $N(x)\cap A=\{v_1, v_2, v_5\}$.

Similarly, $x\in Y$ or $Z$ if $x\sim v_7$.

So, we suppose that $x\not\sim v_2$ and $x\not\sim v_7$. Then, $x\not\sim v_3$ and $x\not\sim v_6$ to forbid a 4-hole on $xv_1v_2v_3x$ or $xv_1v_7v_6x$. Since $|N(x)\cap A|\geq 2$, we have that $x\sim v_4$ or $x\sim v_5$. If $x\sim v_4$ and $x\sim v_5$, then $x\in Z$. Otherwise, $x\in X$. Therefore, (M1) holds.

\medskip

\noindent{\bf (M2)} Both $X_i\cup Z_i$ and $Y_i$ are cliques.

Let $x$, $x'\in X_i\cup Z_i$. If $x\not\sim x'$, then $xv_ix'v_{i+3}x$ is a 4-hole. Let $y$, $y'\in Y_i$. If $y\not\sim y'$, then $yv_iy'v_{i+2}y$ is a 4-hole. This proves (M2).

\medskip

\noindent{\bf (M3)} $Y_i$ is complete to $Y_{i+1}\cup Y_{i+6}$.

Suppose to its contrary, and let $x\in Y_i$ and $y\in Y_{i+1}\cup Y_{i+6}$. If $x\not\sim y$, then $G[\{x, y, v_i, v_{i+1}, v_{i+2}\}]$ is a gem, which leads to a contradiction. Therefore, (M3) holds.

\medskip

\noindent{\bf (M4)} $Y_i$ is anticomplete to $Y_{i+2}\cup Y_{i+3}\cup Y_{i+4} \cup Y_{i+5}$.

Suppose that (M4) does not hold. By symmetry, let $x\in Y_1$ and $y\in Y_{3}\cup Y_{4}\cup Y_{5} \cup Y_{6}$ such that $x\sim y$. If $y\in Y_{3}$, then $G[\{x, y, v_{2}, v_{3}, v_{4}\}]$ is a gem. If $y\in Y_{6}$, then $G[\{x, y, v_1, v_{2}, v_{7}\}]$ is a gem. If $y\in Y_{4}$, then $xv_3v_{4}yx$ is a 4-hole. If $y\in Y_{5}$, then $xv_1v_7yx$ is a 4-hole. All are contradictions. This proves (M4).

\medskip

\noindent{\bf (M5)} Either $X_i=\emptyset$ or $X_{i+2}\cup X_{i+5}=\emptyset$.

If it is not the case, we may choose, without loss of generality, that $u\in X_1$ and $v\in X_3$ then $uvv_3v_4u$ is a 4-hole if $u\sim v$, and $G[\{v, u, v_1, v_2, v_4, v_5, v_6\}]=P_7$ otherwise, which leads to a contradiction. This proves (M5).

\medskip

\noindent{\bf (M6)} $X_i$ is anticomplete to $X_{i+1}\cup X_{i+3}\cup X_{i+4}\cup X_{i+6}$.

Suppose that (M6) does not hold. By symmetry, let $x\in X_1$ and $y\in X_{2}\cup X_{4}\cup X_{5}\cup X_{7}$ such that $x\sim y$. If $y\in X_{2}$, then $xyv_2v_1x$ is a 4-hole. If $y\in X_{4}$, then $xyv_7v_1x$ is a 4-hole. If $y\in X_{5}$, then $xyv_5v_4x$ is a 4-hole. If $y\in X_7$, then $xyv_3v_4x$ is a 4-hole. All are contradictions. This proves (M6).

\medskip

\noindent{\bf (M7)} $X_i$ is complete to $Y_{i+2}\cup Y_{i+6}$.

If (M7) does not hold, we may choose by symmetry that $x\in X_i$ and $y\in Y_{i+2}\cup Y_{i+6}$ such that $x\not\sim y$, then $G[\{x, y, v_i, v_{i+1}, v_{i+2}, v_{i+4}, v_{i+5}\}]=P_7$ whenever $y\in Y_{i+2}$, and $G[\{x, y, v_{i+1}, v_{i+3}, v_{i+4}, v_{i+5}, v_{i+6}\}]=P_7$ whenever $y\in Y_{i+6}$, which leads to a contradiction. Therefore, (M7) hlds.

\medskip

\noindent{\bf (M8)} $X_i$ is anticomplete to $Y_i\cup Y_{i+1}\cup Y_{i+3}\cup Y_{i+4}\cup Y_{i+5}$.

Suppose that (M8) does not hold. Without loss of generality, let $x\in X_1$ and $y\in Y_1\cup Y_2\cup Y_4\cup Y_5\cup Y_6$ such that $x\sim y$. If $y\in Y_1$, then $G[\{x, y, v_1, v_2, v_3\}]$ is a gem. If $y\in Y_2$, then $xyv_2v_1x$ is a 4-hole. If $y\in Y_4$, then $G[\{x, y, v_4, v_5, v_6\}]$ is a gem. If $y\in Y_5$, then $xyv_5v_4x$ is a 4-hole. If $y\in Y_6$, then $G[\{x, y, v_1, v_6, v_7\}]$ is a gem. All are contradictions. This proves (M8).

\medskip

\noindent{\bf (M9)} Either $Z_i=\emptyset$ or $Z_{i+3}\cup Z_{i+4}=\emptyset$.

If this is not true, we may choose by symmetry that $z_1\in Z_1$ and $z_4\in Z_4$, then $z_1v_4z_4v_1z_1$ is a 4-hole whenever $z_1\not\sim z_4$, and $G[\{z_1, z_4, v_1, v_4, v_5\}]$ is a gem whenever $z_1\sim z_4$. This proves (M9).

\medskip

\noindent{\bf (M10)} $Z_i$ is anticomplete to $Z_{i+1}\cup Z_{i+2}\cup Z_{i+5}\cup Z_{i+6}$.

Suppose that (M10) does not hold. By symmetry, let $x\in Z_1$ and $y\in Z_2\cup Z_3\cup Z_6\cup Z_7$ such that $x\sim y$. If $y\in Z_2\cup Z_6$, then  $xyv_2v_1x$ is a 4-hole. If $y\in Z_3\cup Z_7$, then $xyv_7v_1x$ is a 4-hole. Therefore, (M10) is true.

\medskip

\noindent{\bf (M11)} Either $Z_i=\emptyset$ or $X_i\cup X_{i+4}\cup X_{i+6}=\emptyset$.

Suppose that (M11) does not hold. Without loss of generality, we choose $u\in Z_1$ and $v\in X_1\cup X_5\cup X_7$. If $u\not\sim v$, then $uv_1vv_4u$ is a 4-hole whenever $v\in X_1$, $uv_1vv_5u$ is a 4-hole whenever $v\in X_5$, and $G[\{u, v, v_1, v_2, v_3, v_5, v_6\}]=P_7$ whenever $v\in X_7$. Therefore, $u\sim v$, which leads to either a gem $G[\{u, v, v_1, v_4, v_5\}]$ if $v\in X_1\cup X_5$, or a 4-hole $uvv_3v_4u$ if $v\in X_7$.  This proves (M11).

\medskip

\noindent{\bf (M12)} $Z_i$ is anticomplete to $X_{i+1}\cup X_{i+2}\cup X_{i+3}\cup X_{i+5}$.

Suppose that (M12) does not hold. By symmetry, choose $x\in Z_1$ and $y\in X_2\cup X_3\cup X_4\cup X_6$ such that $x\sim y$. If $y\in X_2$, then $xyv_2v_1x$ is a 4-hole. If $y\in X_3$, then $xyv_3v_4x$ is a 4-hole. If $y\in X_4$, then $xyv_7v_1x$ is a 4-hole. If $y\in X_6$, then $xyv_6v_5x$ is a 4-hole. All are contradictions. This proves (M12).

\medskip

\noindent{\bf (M13)} $Z_i$ is complete to $Y_{i+2}\cup Y_{i+3}\cup Y_{i+6}$.

Without loss of generality, let $z\in Z_1$ and $y\in Y_3\cup Y_4\cup Y_7$. If $z\not\sim y$, then $G[\{z, y, v_3, v_4, v_5\}]$ is a gem whenever $y\in Y_3$, $G[\{z, y, v_4, v_5, v_6\}]$ is a gem whenever $y\in Y_4$, and $G[\{z, y, v_2, v_3, v_5, v_6, v_7\}]=P_7$ whenever $y\in Y_7$. Therefore, $z\sim y$, and thus (M13) holds.

\medskip

\noindent{\bf (M14)} $Z_i$ is anticomplete to $Y_i\cup Y_{i+1}\cup Y_{i+4}\cup Y_{i+5}$.

By symmetry,  choose $z\in Z_1$ and $y\in Y_1\cup Y_2\cup Y_5\cup Y_6$. If  $z\sim y$, then $G[\{z, y, v_1, v_2, v_3\}]$ is a gem whenever $y\in Y_1$, $zyv_2v_1z$ is a 4-hole whenever $y\in Y_2$, $G[\{z, y, v_4, v_5, v_6\}]$ is a gem whenever $y\in Y_5$, and $G[\{z, y, v_1, v_6, v_7\}]$ is a gem whenever $y\in Y_6$. Therefore, $z\not\sim y$, and thus (M14) holds.

\medskip

Let $Y'_i=Y_{i-1}\cup \{v_i\}$ for each $i\in\{1, \cdots, 7\}$ and $Y'=Y'_1\cup\cdots\cup Y'_7$. By (M2), (M3) and (M4), we have that $G[Y]$ is a clique blowup of $C_7$ and hence $G[Y']$ is a nonempty clique blowup of $C_7$. Follow from (M1) to (M14), $G[N(A)\cup A]$ is a clique blowup of some graph like the one as shown in Figure 4.

\begin{figure}[htbp]\label{fig-4}
	\begin{center}
		\includegraphics[width=3.5cm]{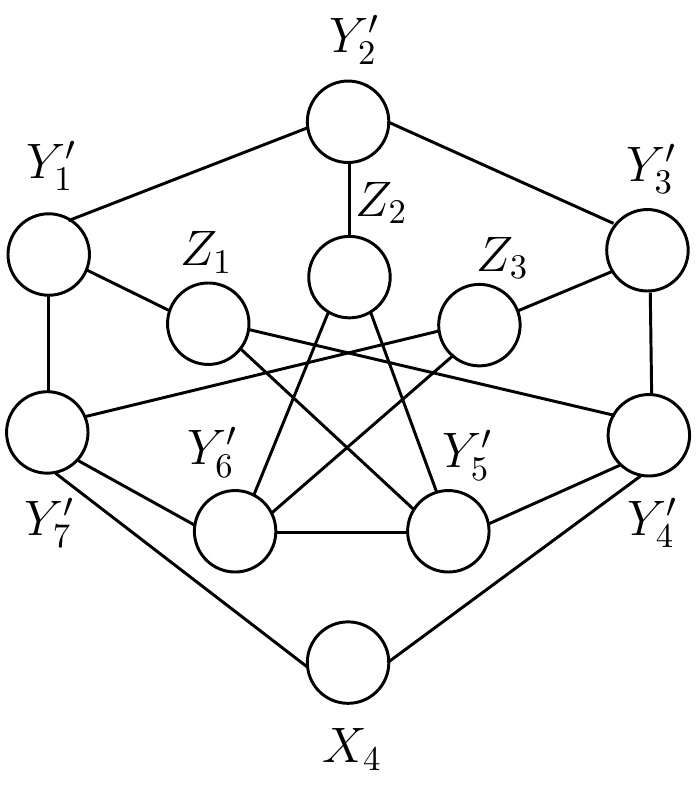}
	\end{center}
	\vskip -25pt
	\caption{One situation of $G[ N(A)\cup A] $.}
\end{figure}

Now, we can prove Theorem~\ref{gem}, by considering  whether $X_i$ is empty or not, to find a bisimplicial vertex of $G$ in $A$.

By (M5), we have that $\{X_1,\cdots,X_7\}$ has at most three nonempty elements. So, there must exist an integer $i\in\{1,\cdots,7\}$ such that $v_i$ is anticomplete to $X$.

Suppose that $Z=\emptyset$. Then, $N(v_i)=(Y_{i-2}\cup\{v_{i-1}\}\cup Y_{i-1})\cup(Y_{i}\cup\{v_{i+1}\})$ by (M1).  Notice that $v_{i-1}$ is complete to $Y_{i-2}\cup Y_{i-1}$, and  $v_{i+1}$ is complete to $Y_i$, and $Y_{i-2}$ is complete to $Y_{i-1}$ by (M3). We have that $Y_{i-2}\cup \{v_{i-1}\}\cup Y_{i-1}$ and $Y_{i}\cup\{v_{i+1}\}$ are two nonempty cliques by (M2), and thus $v_i$ is a bisimplicial vertex of $G$.

Therefore, we suppose that $Z\ne\emptyset$.
Without loss of generality, suppose $Z_1\ne\emptyset$. By (M9) and (M11), $Z_4\cup Z_5\cup X_1\cup X_5\cup X_7=\emptyset$, and so
$$X=X_2\cup X_3\cup X_4\cup X_6, \mbox{ and } Z=Z_1\cup Z_2\cup Z_3\cup Z_6\cup Z_7.$$

Firstly, we prove that

\begin{equation}\label{gem-1}
	\mbox{if $X_2=\emptyset$, then $v_5$ is a bisimplicial vertex of $G$.}	
\end{equation}

Suppose that $X_2=\emptyset$. By (M1), we can deduce that $N(v_5)=(Y_4\cup Y_3\cup Z_1\cup \{v_4\})\cup (Y_5\cup Z_2\cup \{v_6\})$. By (M3) and (M13), $Y_3$ is complete to $Y_4$, $Z_1$ is complete to $Y_3\cup Y_4$, and  $Y_5$ is complete to $Z_2$. By the definition of $Z_i$ and $Y_i$, $v_4$ is complete to $Y_4\cup Y_3\cup Z_1$, and $v_6$ is complete to $Y_5\cup Z_2$. It follows from (M2) that  $Y_4\cup Y_3\cup Z_1\cup \{v_4\}$ and $Y_5\cup Z_2\cup \{v_6\}$ are two nonempty cliques. Therefore, $v_5$ is a bisimplicial vertex of $G$. This proves (\ref{gem-1}).

Next, we prove that

\begin{equation}\label{gem-2}
	\mbox{if $X_4=\emptyset$, then $v_4$ is  a bisimplicial vertex of $G$.}	
\end{equation}

Suppose that $X_4=\emptyset$.  Then, $N(v_4)=(Y_3\cup Y_2\cup Z_7\cup \{v_3\})\cup (Y_4\cup Z_1\cup \{v_5\})$ by (M1). By (M3) and (M13), $Y_2$ is complete to $Y_3$, $Z_7$ is complete to $Y_2\cup Y_3$, and $Z_1$ is complete to $Y_4$. By (M2) and by the definition of $Z_i$ and $Y_i$, we have that $Z_1$, $Z_7$, $Y_2$, $Y_3$ and $Y_4$ are all cliques, and $v_3$ is complete to $Y_3\cup Y_2\cup Z_7$, and $v_5$ is complete to $Y_4\cup Z_1$. Therefore, $Y_3\cup Y_2\cup Z_7\cup \{v_3\}$ and $Y_4\cup Z_1\cup\{v_5\}$ are two cliques, and so $v_4$ is a bisimplicial vertex of $G$. This proves (\ref{gem-2}).

By (M5), we have that one of $X_2$ and $X_4$ is empty. It follows immediately from (\ref{gem-1}) and (\ref{gem-2}) that $G$ must have a bisimplicial vertex. This completes the proof of Theorem~\ref{gem}. \qed

\medskip


It is certain that $G_2$ defined in Section 2 is ($P_7$, gem)-free and contains a 7-hole, and it is easy to check that $G_2$ does not satisfy the conclusion of Theorem~\ref{gem}.

Let $A_3=z_1z_2z_3z_4z_5z_6z_7z_1$ be a 7-hole. We use $G_5$ to denote the graph obtained from $A_3$ by adding seven vertices $a_1$, $a_2$, $\cdots$, $a_7$ such that $N(a_i)\cap (V(A_3)\cup\{a_1, a_2, \cdots, a_7\}) =\{z_i, z_{i+3}, z_{i+4}, a_{i+3}, a_{i+4}\}$ for each $i\in\{1, \cdots, 7\}$ (see Figure 5). Obviously, $G_5$ is $(P_7, C_4$)-free and contains a 7-hole and thus the clique blowup of $G_5$ is $(P_7, C_4$)-free and contains a 7-hole. Moreover, one can easily check that the clique blowup of $G_5$ does not satisfy the conclusion of Theorem~\ref{gem}.

\begin{figure}[htbp]\label{fig-5}
	\begin{center}
		\includegraphics[width=4cm]{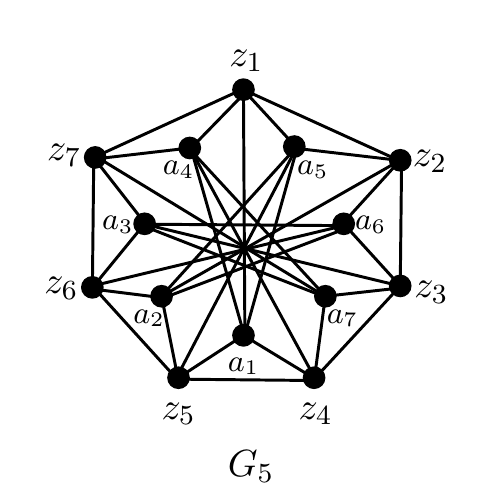}
	\end{center}
	\vskip -25pt
	\caption{Illustration of $G_5$.}
\end{figure}

We use $G_6$ to denote the $t$-size  clique blowup of $G_3$, defined in Section 2, with $t\geq 1$. It is certain that $G_6$ is ($C_4$, gem)-free and contains a 7-hole, and we can check that $G_6$ does not satisfy the conclusion of Theorem~\ref{gem}.

These graphs show that the requirements $P_7$-free, $C_4$-free, and gem-free in Theorem~\ref{gem} are all necessary.

\medskip

Now, we can prove Corollary~\ref{gem*} by induction on $|V(G)|$. Let $G$ be a $(P_7, C_4$, gem)-free graph. We may assume that $G$ is connected and has no clique cutsets. By Theorem~\ref{gem}, we have that $G$ has a bisimplicial vertex or $G$ is the clique blowup of the Petersen graph. By Lemma~\ref{Petersen}, it follows immediately that $\chi(G)\le2\omega(G)-1$.

\bigskip

\noindent{\bf Remarks}. Recall that we use $P_k$ to denote a path on $k$ vertices, and use diamond to denote the graph obtained from a $K_4$ by removing an edge. In 1998, Randerath \cite{BR1998} proved that $\chi(G)\le \omega(G)+1$ for every ($P_5$, diamond)-free graph.  Cameron, Huang and Merkel \cite{CHM18} proved that $\chi(G)\leq \omega(G)+3$. for every ($P_6$, diamond)-free graph. Schiermeyer \cite{SR2023} asked a question: Is it true that there exists a constant $c$ such that $\chi(G)\le \omega(G)+c$ for every $(P_7$, diamond)-free graph$?$ This is still open.

\bigskip

\noindent{\bf Acknowledgement}: We thank T. Karthick for informing us reference \cite{SR2023}, and thank I. Schiermeyer for helpful suggestions.


\begin{thebibliography}{9999}

\bibitem{BM08} J. A. Bondy, U. S. R. Murty, Graph Theory, Springer, New York, 2008.

\bibitem{CHM18} K. Cameron, S. Huang, O. Merkel, An optimal $\chi$-bound for ($P_6$, diamond)-free graphs, J. of Graph Theory, 97 (2021) 451-465.

\bibitem{CHPS20} K. Cameron, S. Huang, I. Penev, V. Sivaraman, The class of $(P_7, C_4, C_5)$-free graphs: Decomposition, algorithms, and $\chi$-boundedness, J. of Graph Theory, 93 (2020) 503-552.

\bibitem{CKB21} S. A. Choudum, T. Karthick, M. M. Belavadi, Structural domination and coloring of some $(P_7, C_7)$-free graphs, Disc. Math., 344 (2021) 112244.

\bibitem{CKS07} S. A. Choudum, T. Karthick, M. A. Shalu, Perfect coloring and linearly $\chi$-bound $P_6$-free graphs, J. of Graph Theory, 54 (2007) 293-306.

\bibitem{CS19} M. Chudnovsky, V. Sivaraman, Perfect divisibility and 2-divisibility, J. of Graph Theory, 90 (2019) 54-60.

\bibitem{FHH18} D. J. Fraser, A. M. Hamel, C. T. Ho\'{a}ng, On the structure of (even-hole, kite)-free graphs, Graphs and Combinatorics, 34 (2018) 989-999.

\bibitem{GSH17} S. Gaspers, S. Huang, Linearly $\chi $-Bounding $(P_6, C_4)$ -Free Graphs, International Workshop on Graph-Theoretic Concepts in Computer Science, Springer, Cham, 2017.

\bibitem{G22} M. Gei{\ss}er, Colourings of $P_5$-free graphs, PhD thesis, 2022.

\bibitem{GHM03} S. Gravier, C. T. Ho\'{a}ng, F. Maffray, Coloring the hypergraph of maximal cliques of a graph with no long path, Disc. Math., 272 (2003) 285-290.

\bibitem{G75} A. Gy\'{a}rf\'{a}s, On Ramsey covering-numbers, Infinite and Finite Sets, 2 (1975) 801-816.

\bibitem{H22} S. Huang, The optimal $\chi $-bound for $(P_7, C_4, C_5)$-free graphs, arXiv:2212.05239, 2022.

\bibitem{KM19} T. Karthick, F. Maffray, Square-free graphs with no six-vertex induced path, SIAM J. on Disc. Math., 33 (2019) 874-909.

\bibitem{LZL22} K. Lan, Y. Zhou, F. Liu, The chromatic number of $(P_6, C_4$, diamond)-free graphs, Bulletin of the Australian Mathematical Society, (2022) 1-10.

\bibitem{M21} S. Mishra, On graphs with no induced bull and no induced diamond, arXiv:2107.03750, 2021.

\bibitem{BR1998} B. Randerath, The Vizing bound for the chromatic number based on forbidden
Pairs (Ph. D. Thesis), Shaker Verlag, Aachen, (1998).

\bibitem{RS2004} B. Randerath, I. Schiermeyer, Vertex colouring and forbidden subgraphs-A survey, Graphs and Combinatorics, 20 (2004) 1-40.

\bibitem{S16} I. Schiermeyer, Chromatic number of $P_5$-free graphs: Reed's conjecture, Disc. Math., 339 (2016) 1940-1943.

\bibitem{SR2019} I. Schiermeyer, B. Randerath, Polynomial $\chi$-binding functions and forbidden induced subgraphs: A survey, Graphs and Combinatorics, 35  (2019) 1-31.

\bibitem{SR2023} I. Schiermeyer, Polynomial $\chi$-binding functions for $P_5$-free graphs, personal communication.

\bibitem{SS18} A. Scott, P. Seymour, A survey of $\chi$-boundedness, J. of Graph Theory, 95 (2020) 473-504.

\bibitem{S74} D. Seinsche, On a property of the class of $n$-colorable graphs, J. of Combinatorial Theory, Ser. B, 16 (1974) 191-193.






\end{thebibliography}
\end{document}